\documentclass[12pt, a4paper]{amsart}
\usepackage{graphics,color,hyperref,verbatim}
\usepackage{amsmath}
\usepackage{amsfonts}
\usepackage{amssymb}
\usepackage{amsthm}
\usepackage{graphicx}
\usepackage{psfrag}
\usepackage{marvosym}
\usepackage{mathrsfs}
\usepackage{enumitem}

\usepackage[a4paper, includeheadfoot, margin=2.2 cm, top=1.2 cm, footskip=1.2 cm]{geometry} 

\usepackage[usenames,dvipsnames]{xcolor} 
\usepackage{tikz, tikz-3dplot, pgfplots}
\usetikzlibrary{decorations.pathreplacing,
                calligraphy,
                matrix}
\usepackage{tikz-cd}
\usetikzlibrary[positioning,patterns] 
\usetikzlibrary{matrix,arrows,decorations.pathmorphing}

\pgfplotsset{compat=1.16}

\usepackage{ytableau}

\usepackage{calligra}
\usepackage{mathrsfs}
\DeclareMathAlphabet{\mathcalligra}{T1}{calligra}{m}{n}
\DeclareFontShape{T1}{calligra}{m}{n}{<->s*[1.5]callig15}{}

\newtheorem{theorem}{Theorem}[section]

\newtheorem{lemma}[theorem]{Lemma}
\newtheorem{proposition}[theorem]{Proposition}
\newtheorem{corollary}[theorem]{Corollary}

\newtheorem{question}[theorem]{Question}

\theoremstyle{definition}
\newtheorem{definition}[theorem]{Definition}

\newtheorem{example}[theorem]{Example}

\newtheorem{remark}[theorem]{Remark}
\newtheorem{theorem-definition}[theorem]{Theorem-Definition}

\numberwithin{equation}{section}

\renewcommand{\AA} {\mathbb{A}}

\newcommand{\DD} {\mathbb{D}}

\newcommand{\MM} {\mathbb{M}}

\newcommand{\PP} {\mathbb{P}}

\newcommand{\ZZ} {\mathbb{Z}}

\def\be{\mathbf{e}}

\newcommand {\shO} {\mathcal{O}}

\newcommand {\sI} {\mathscr{I}}

\newcommand {\sO} {\mathscr{O}}

\newcommand {\fom}  {\mathfrak{m}}

\newcommand{\blank}{\underline{\hphantom{A}}}

\newcommand {\codim} {\operatorname{codim}}

\newcommand {\Ext} {\operatorname{Ext}}
\newcommand{\sExt}{\mathscr{E} \kern -1pt xt}

\newcommand {\GL} {\operatorname{GL}}

\newcommand {\Gr} {\mathrm{Gr}}

\newcommand{\Hilb}{\mathrm{Hilb}}

\newcommand {\Hom} {\operatorname{Hom}}
\newcommand {\sHom}{\mathscr{H}\kern-5pt\mathcalligra{om}}

\newcommand {\kk} {\Bbbk}

\newcommand {\Mat} {\operatorname{Mat}}

\newcommand {\rank} {\operatorname{rank}}

\newcommand{\sTor}{\mathscr{T} \kern -3pt or}

\newcommand {\Bl} {\operatorname{Bl}}

\newcommand{\BN}{\operatorname{BN}}

\newcommand{\Grass}{\operatorname{Grass}}

\def\onto{\ensuremath{\twoheadrightarrow}}

\newcommand\scalemath[2]{\scalebox{#1}{\mbox{\ensuremath{\displaystyle #2}}}}

\title[BN for $\Hilb_n$]{Brill--Noether theory of Hilbert schemes of points on surfaces}

\author[A. \ BAYER]{Arend Bayer}
\address{School of Mathematics, University of Edinburgh, JCMB, Peter Guthrie Tait Road, Edinburgh EH9 3FD, UK.} 
\email{arend.bayer@ed.ac.uk} 

\author[H.C. \ CHEN]{Huachen Chen}
\address{Department of Mathematics,  University of California Santa Barbara, CA 93106, USA.} 
\email{hchen@math.ucsb.edu}

\author[Q.Y.\ JIANG]{Qingyuan Jiang}
\address{School of Mathematics, University of Edinburgh, JCMB, Peter Guthrie Tait Road, Edinburgh EH9 3FD, UK.} 
\curraddr{Department of Mathematics,
The Hong Kong University of Science and Technology, Clearwater Bay, Kowloon, Hong Kong.} 
\email{jiangqy@ust.hk}

\begin{document}

\begin{abstract} We show that Brill--Noether loci in Hilbert scheme of points on a smooth connected surface $S$ are non-empty whenever their expected dimension is positive, and that they are irreducible and have expected dimensions. More precisely, we consider the loci of pairs $(I, s)$ where $I$ is an ideal that locally at the point $s$ of $S$ needs a given number of generators. 

We give two proofs. The first uses Iarrobino's description \cite{Ia} of the Hilbert--Samuel stratification of local punctual Hilbert schemes, and the second  is based on induction via birational relationships between different Brill--Noether loci given by nested Hilbert schemes.
 \end{abstract}

\maketitle
 
 \section{Introduction}
Given any pair $M, N$ of moduli spaces of sheaves or complexes on a given variety, one may consider Brill--Noether loci in $M \times N$ determined by the dimension of the space of morphisms between the corresponding objects.
\begin{question} \label{question}
In which situations are Brill-Noether loci in $M \times N$ of expected dimension?
\end{question}
This question was, in a sense, first considered in Lazarsfeld's proof of Brill--Noether for curves on K3 surfaces \cite{Lazarsfeld:BN}, and many variants of classical Brill--Noether for curves on surfaces can be formulated in this way; see also \cite{B, BL17}. A dual version is a crucial ingredient for Le Potier's Strange duality \cite{MO:strange}. In general, such Brill-Noether loci can be seen as generalisations of nested Hilbert schemes, which highlights a wide range for potential implications, by giving relations between the birational geometry, the cohomology, Chow groups, motives or the derived categories of different moduli spaces.
 
In this paper, we consider the fundamental case where $M$ is the Hilbert scheme of points on a surface, and $N$ is the surface itself.
 Let $S$ be a smooth irreducible surface over a field $\kk$, and $\Hilb_n(S)$ the Hilbert scheme of ideals $I \subset \sO_S$ of colength $n$. For each integer $r \ge 0$, the {\em $r$th Brill--Noether locus} is the locus in $\Hilb_n(S) \times S$ of pairs $(I,p)$ with $\dim_{\kappa(p)} \Hom(I, \kappa(p)) \geq r+1$; equivalently, the locus where the minimal number of generators of $I$ locally at $p$ is at least $r+1$, which by Nakayama's Lemma is
	\[ \BN_{r,n} : = \{(I,p) \mid \dim_{\kappa(p)} (I \otimes_{\sO_S} \kappa(p)) \ge r+1\} \subset  \Hilb_n(S) \times S. \]
Then $\BN_{0,n} = \Hilb_n(S) \times S$, and $\BN_{1,n} = Z_n := \{(I,p) \mid p \in V(I)\}\subset \Hilb_n(S) \times S$ is the universal subscheme.

The main result of our paper completely answers Question \ref{question} in our setting:
\begin{theorem} \label{thm:main} 
For each $n \ge 1$ and $r \ge 0$, the Brill--Noether locus $\BN_{r,n}$ is Cohen--Macaulay, irreducible and of expected codimension
	$$\codim (\BN_{r,n}, \, \Hilb_n(S) \times S) = r(r+1).$$ 
Equivalently, $\dim \BN_{r,n} = \rho_{r,n} := 2n+2 - r(r+1).$
Moreover, $\BN_{r,n} \ne \emptyset$ iff $ \rho_{r,n} \ge 2$. 
\end{theorem}

\begin{remark}
One can also observe that $\BN_{r,n}$ is the closure of the locus of pairs of the form $(\fom_p^r \cdot J, p)$ where $p \in S$, $\fom_p$ is the maximal ideal of $p$, and $J$ is an ideal of $\frac{\rho_{r,n}}2-1$ points supported away from $p$. Indeed, this clearly gives a locus of dimension $\rho_{r, n}$, and thus the claim follows by the irreducibility of $\BN_{r,n}$. 
\end{remark}

\subsection{Previous results}
 Ellingsrud and Str{\o}mme \cite[Proposition 2.2]{ES}  proved the bound:
	$$\codim (\BN_{r,n}, \, \Hilb_n(S) \times S) \ge 2r \quad \text{for all} \quad r \ge 0.$$
It was used repeatedly to study the geometry of nested Hilbert schemes: in \cite{ES} to prove irreducibility, in \cite{RY} to study their nef cones, and in \cite{JL18,J19}
to study their derived categories and Chow groups. The bound was improved by Ryan and  Taylor \cite{RT}:
    $$\codim (\BN_{r,n}, \, \Hilb_n(S) \times S) \ge \binom{r+1}{2} +1 \quad \text{for all} \quad r \ge 1, n \ge  \binom{r+1}{2}.$$
    
Finally, if $I$ has $r+1$ generators at $p$ then the \emph{socle} of $\shO_{V(I), p}$ is $r$-dimensional (see \cite[Lemma 2.1]{Song}). Thus the non-emptiness statement of Theorem~\ref{thm:main} is equivalent to the sharp bound on the dimension of the socle established in \cite[Theorem 1.2]{Song}. One direction of this bound,  the emptiness of $\BN_{r,n}$ for $\rho_{r, n}\le 0$, was first established by Haiman \cite[Proposition 3.5.3]{Hai2} in his analysis of the singularities of the isospectral Hilbert scheme.

\subsection{Proofs}
We present two proofs. The first, in Section \ref{sec:local}  is based on explicit resolutions of $I$  over Hilbert--Samuel strata of local punctual Hilbert schemes given by  \cite{Bri, Ia}.

The second proof, in Section \ref{sec:biration} is an inductive argument based on birational relations among the various $\BN_{r,n}$ induced by nested Hilbert schemes. We expect that this method will be useful in answering Question \ref{question} in much bigger generality.

Our method also gives the irreducibility of one type of nested Hilbert scheme, see Remark~\ref{rem:nestedHilbert}. In general, irreducibility of nested Hilbert schemes is an open problem; see 
\cite{GRS}, \cite{RT} and \cite[Section 2.A]{Add} for recent results on this topic.

\subsection*{Acknowledgement} 
We thank the referees for their many helpful comments and suggestions.
A.B.~ and Q.J.~were supported by EPSRC  grant EP/R034826/1, and by the ERC Grant ERC-2018-CoG-819864-WallCrossAG. 
Q.J. is also funded by the Deutsche Forschungsgemeinschaft (DFG, German Research Foundation) under Germany's Excellence Strategy – EXC-2047/1 – 390685813.

\section{Brill--Noether loci via local Hilbert schemes}\label{sec:local}

In this section, we study Brill-Noether loci of local Hilbert schemes (by which we mean that the entire subscheme is supported at a given point). Our proof is based on explicit  coordinate charts of their Hilbert--Samuel strata constructed by Iarrobino \cite{Ia} and Brian\c{c}on \cite{Bri}. This leads to our first proof of Theorem \ref{thm:main}.

Since it is enough to prove Theorem~\ref{thm:main} after base change to the algebraic closure, we will from now on assume that $\kk$ is algebraically closed.

\subsection{Hilbert--Samuel stratification}
We follow the terminology and convention of \cite{Ia}. 
Let $A=\kk[[x,y]]$ be the ring of power series in two variables $x$ and $y$, where $\kk$ is an algebraically closed field, and let $\fom = (x,y)$ denote the maximal ideal. For any ideal $I \subset A$, the natural grading of $A$ induces the Hilbert--Samuel function of $A/I$ given by 
	$$\chi_{A/I}(i) = \dim_{\kk}\left(\frac{A}{I+\fom^{i+1}}\right) = \dim_{\kk}\left(\frac{A/I}{\fom^{i+1} (A/I)} \right) \qquad i \in \ZZ_{\ge 0}.$$
It is usually convenient to consider the step function of the Hilbert--Samuel function, as follows.

 \begin{definition} \label{def:type}
     Given an ideal $I \subset A$, its \emph{type} $T(I)  =(t_0, t_1, t_2, \dots )$ is the sequence
  \[
   t_j 
   = \dim_\kk \left( \frac{I + \fom^j}{I + \fom^{j+1}} \right), 
  \]
  and its \emph{order} $d$ is determined by 
  \[ d(I) = \sup\{k \in \ZZ_{\ge 0} \mid I \subset \fom^k\}. \]
 \end{definition}
For $i \ge 1$, the short exact sequence
\[ 0 \to \frac{I + \fom^{i}}{I + \fom^{i+1}} \to  \frac{A}{I + \fom^{i+1}} \to  \frac{A}{I + \fom^{i}} \to 0 \]
    shows that  $t_i(I) = \chi_{A/I}(i) - \chi_{A/I}(i-1)$. Moreover, if $I$ has colength $n$ and order $d$, then its type satisfies
\begin{eqnarray} \label{eq:typeproperties}
T = (1, 2, 3, \dots, d, t_{d}, t_{d+1}, \dots, 0, 0, 0, \dots) \quad &\text{where}& \quad
d \ge t_d \ge t_{d+1} \ge \cdots \ge 0 \nonumber \\
\quad  &\text{and}& \quad |T| = \sum_j t_j = n. 
\end{eqnarray} 
Conversely, for any type satisfying \eqref{eq:typeproperties}, there exists an ideal of type $T$, see Example \ref{ex:normalideal}. 

\begin{example} \label{ex:grassmannstratum}
    Let $d$ be such that $0 \le \ell:= n - \frac{d(d+1)}2 \le d$. Then $I$ belongs to the Grassmannian stratum $\fom^{d+1} \subset I \subset \fom^{d}$ parametrised by $(d+1-\ell)$-dimensional subspaces of $\fom^d/\fom^{d+1}$ if and only if  $I$ is of type $(1, 2, 3, \dots, d, \ell, 0, \dots)$.
\end{example}

It is often convenient to encode the type via the following data.
\begin{definition} \label{def:ei}
Given a type $T$ of order $d$, we define the {\em jumping indices of $T$} by 
	\begin{equation*} 
	e_j = 
	\begin{cases}
	t_{j-1} - t_j  & \text{~if~} j\ge d \\
	 0 &  \text{~if otherwise}. 
	\end{cases}
	\end{equation*}
\end{definition}
Then $e_j \ge 0$ for each $j$, $e_j=0$ for $j  \ge n+1$ or $j < d$, and $\sum e_j = d$. Clearly, $d$ and $T$ are determined by the jumping indices.

\begin{theorem}[Hilbert--Samuel Stratifications; \cite{Bri, Ia}] \label{thm:stratification}
Let $\kk$ be an algebraically closed field, let $A=\kk[[x,y]]$, $n \ge 2$, and let $\Hilb_n(A) $ denote the local punctual Hilbert scheme, with reduced scheme structure. For each type $T$ with $|T|=n$, let $Z_T$ denote the subset of $\Hilb_n(A)$ consisting of ideals $I$ of type $T$, and let $e_j$ denote the jumping indices of $T$.

\begin{enumerate}
	\item There is a decomposition of $\Hilb_n(A)$ into a disjoint union 
	$$\Hilb_n(A) = \bigsqcup_{|T| = n} Z_T,$$
where $T$ runs through all types with $|T|=n$ satisfying \eqref{eq:typeproperties}.
	\item For each type $T$ satisfying \eqref{eq:typeproperties}, the stratum $Z_T$ is a locally closed subset of $\Hilb_n(A)$, which is nonempty, smooth, rational, connected, of dimension 
		$$
		\dim Z_T = n - \sum_{j \ge d} \frac{e_j(e_j+1)}{2} = n - d - \sum_{j \ge d} \frac{e_j (e_j-1)}{2}.
		$$
\end{enumerate}
\end{theorem}
\begin{proof}
 Since the Hilbert--Samuel function is upper-semicontinuous on $\Hilb_n(A)$, the difference function is constructible, and  thus each $Z_T$ is locally closed. 

 By \cite[Theorem 3.13]{Ia}, each $Z_T$ is irreducible, rational and nonsingular. The claim about the dimension of $Z_T$ follows from \cite[Theorem 2.12]{Ia} (see also \cite[Theorem III.3.1]{Bri}).
\end{proof}

Despite being commonly referred to as the Hilbert--Samuel stratification, the decomposition of Theorem \ref{thm:stratification} does {\em not} always satisfy the condition that the closure of a stratum is a union of strata.

\begin{example}[Curvilinear Strata]
For each $n \ge 2$, there is a unique type of order $d=1$:
	$$T_{n, \rm curv} = (\underbrace{1,1,1, \ldots, 1,1}_{\text{$n$ terms}}).$$
The corresponding stratum $Z_{n, \rm curv} := Z_{T_{n, \rm curv}} \subset \Hilb_n(A)$ is  the {\em curve-linear stratum} of $\Hilb_n(A)$: an ideal has type $T_{n, \rm curv}$ if and only if it has order 1, i.e. if the associated zero-dimensional subscheme is contained in the germ of a smooth curve defined by $f \in I$, $f \notin \fom^2$. Moreover, $I$ is determined by $f$ via $I = f + \fom^n$.

Using affine coordinate charts, it is easy to see that $Z_{n, \rm curv}$ is smooth, connected and has dimension $n-1$.
The curvilinear stratum $Z_{n, \rm curv} \subset \Hilb_n(A)$ is open and dense (Remark \ref{rmk:irreducibility}), and each ideal $I \in Z_{n, \rm curv}$ can be generated by two elements.
\end{example}

\begin{remark} \label{rmk:irreducibility} 
One primary application  of Theorem \ref{thm:stratification} in \cite{Bri, Ia} was to establish the irreducibility of punctual Hilbert schemes. Concretely, using a deformation argument, Briancon \cite[Theorem V.3.2; Corollary V.3.3]{Bri} and Iarrobino  \cite[\S 5]{Ia}
showed that $\overline{Z_{n, \rm curv}} = \Hilb_n(A)$. Consequently, $\Hilb_n(A)$ is irreducible of dimension $n-1$.
\end{remark}

\subsection{Normal patterns and affine charts of Hilbert--Samuel strata}
\label{sec:normal}
The notion of normal patterns gives rise to affine covers of each stratum $Z_T$ associated with a type $T$. Concretely, given a type $T = (t_j)$, the {\em normal pattern} $P$ of type $T$
is the set of monomials:
	$$P = \bigcup_{j \ge 0} P_j \qquad P_j = \{x^{j-t} y^t \mid 0 \le t \le t_j - 1\}.$$
The normal pattern $P$ associated with a type $T$ can be visualised as Young diagram $\Delta(P)$ such that $(i,j) \in \Delta(P) \iff x^i y^j \in P$. For example, in the case where $T = (1,2,3,2,2,0,0, \ldots)$, we can depict $P$ as 
 	\[
P= 
	\begin{ytableau}
y^2  \\
y & xy & x^2y & x^3y\\
1 & x & x^2 & x^3 & x^4
\end{ytableau}
\]
\begin{remark} \label{rmk:typevski}
The Young diagram $\Delta(P)$ is determined by the property that the row lengths $k_0 > k_1 > \dots > k_d = 0$ give a strictly decreasing partition of $n$, and that it contains $t_j$ monomials of degree $j$. 

Conversely, given any Young diagram $\Delta(P)$ with strictly decreasing row lengths, $t_j$ is the number of monomials in $P$ of degree $j$, whereas the jumping indices $e_j$ are determined as follows: $e_j=0$ if $j < d$ or $j >k_0$, and if $d \le j \le k_0$, $e_j$ is the number of degree $j$ monomials in the sequence 
    \[x^{k_0} y^0, x^{k_1} y^1, \dots, x^{k_{d-1}}y^{d-1}.\]
\end{remark}

\begin{definition}[Affine Charts Associated with Normal Patterns]
\label{def:ZP}
Let $P=P(x,y)$ with type $T=T(P)$ be a normal pattern. We let $Z_P \subset Z_T$  denote the subset of ideals $I$ satisfying the following equivalent conditions (\cite[Lemma 1.4]{Ia}):
\begin{enumerate}
			\item For all $j$, $\langle P \cap \fom^j \rangle \oplus (I \cap \fom^j) = \fom^j$.
			\item $\langle P \rangle \cap I = 0$ and $T(P) = T(I)$, where $\langle P \rangle$ denotes the $\kk$-linear span of $P$.
		\end{enumerate}
\end{definition}
Then $Z_P$ is a Zariski open subscheme of $Z_T$, isomorphic to an affine space (\cite[Propositions 2.5 \& 2.8]{Ia}). The proof of Theorem \ref{thm:main} will use an explicit parametrisation of $Z_P$.

\begin{example} \label{ex:normalideal}
    Let $k_0 > k_1 > k_2 > \dots> k_{d-1}$ be the partition of $n$ corresponding to $\Delta(P)$ for a type $T$ and order $d$, and set $k_d = 0$. Let
    $u_s = x^{k_s} y^s$ for $0 \leq s \leq d$. Then the monomial ideal $I = (u_0, \dots, u_d)$ is contained in $Z_P$. 
\end{example}

In characteristic zero or large characteristics, when we vary the system of parameters $(x,y)$ linearly, the affine spaces $Z_{P}$ of Definition \ref{def:ZP} form an open covering of the stratum $Z_T$:
\begin{proposition}[{\cite[Proposition 3.2 \& Corollary 3.3]{Ia}}]
\label{prop:affinecover}
Let $T$ be a type with $|T|=n$. Assume that either ${\rm char}(\kk) =0$ or ${\rm char}(\kk) \ge |T|=n$. 
Then $Z_T$ is the union of a finite number of translates of $Z_P$ under the action of $\mathrm{GL}(2, \kk)$.
\end{proposition}

The above statements are no longer valid in low characteristic cases.
This proposition is the only place where we need the characteristic assumption in our first proof of the main theorem.

The monomial ideal of Example \ref{ex:normalideal} has a resolution of the form
\begin{equation} \label{eq:normalresolution}
A^d \xrightarrow[]{\MM_P \cdot } A^{d+1} \xrightarrow[]{(u_0, \dots, u_d) \cdot } I
\end{equation}
where the $(d+1) \times d$-matrix $\MM_P$ encodes the obvious relations $y \cdot u_{s-1} = x^{k_{s-1} - k_{s}} \cdot u_s$:
\[
\left(\MM_P\right)_{ii} = -y, \quad \left(\MM_P\right)_{(j+1)j} = x^{k_{j-1} - k_{j}}, \quad
\left(\MM_P\right)_{ij} = 0 \text{ if $j \neq i, i-1$.}
\]
By Nakayama's Lemma, as $\MM_P(0, 0) = 0$, this shows that $I$ needs $d+1$ generators. Note also that, up to signs, the $u_i$ are precisely the $d \times d$-minors of $\MM_P$ .

\subsection{Affine parametrization}
\label{sec:2ndparameters}
We will use a parametrisation of $Z_P$ due to Iarrobino that is obtained by deforming $I$ via deforming the matrix $\MM_P$ and thus the resolution \eqref{eq:normalresolution}.

\begin{proposition}[{\cite[Proposition 4.17]{Ia}}\footnote{Notice that the inequality ``$\mu \le w_{\max\{i,j\}}$" in the formula \cite[(4.14)]{Ia} should be ``$\mu < w_{\max\{i,j\}}$".}] 
\label{prop:Iaparameters}
For a normal pattern $P$ of type $T$, where $T$ has order $d$ and $|T|=n$, we consider
$(d+1)\times d$ matrices $\beta$ with entries in $\kk[x]$ satisfying the following constraints:
\begin{enumerate}
    \item  \label{enum:betavanish} $\beta_{ij} = 0$ if $i > j$;
    \item \label{enum:betadegree}  for $i \leq j$, the entry $\beta_{ij}$ is a polynomial of degree $k_{j-1} - k_j - 1$; and
    \item \label{enum:beta0} $\beta_{ij}(0) = 0$ if $k_{j-1} + j = k_{i-1} + i$.
\end{enumerate}

Let $I(\beta)$ be the ideal generated by the $d \times d$-minors of $\MM_P + \beta$. Then $I(\beta) \in Z_P$, and 
\begin{equation} \label{eq:Ibetaresolution}
  A^d \xrightarrow{\MM_P + \beta} A^{d+1} \to I(\beta)  
\end{equation}
is a resolution for $I(\beta)$. Conversely, any ideal in $Z_P$ is of the form $I(\beta)$ for a unique matrix $\beta$ satisfying the conditions (\ref{enum:betavanish})--(\ref{enum:beta0}) above.
\end{proposition}

We first consider condition (\ref{enum:beta0}) in more detail. It applies exactly when $i, j$ belong to the same group $i', \dots, i' + e_{i'+k_{i'}}-1$, for some $i'$, occurring in Remark~\ref{rmk:typevski}. In other words, $\beta(0)$ is a matrix whose bottom row vanishies, and whose top $d \times d$-block satisfies the following condition.

\begin{definition}
\label{def:Mat_e}
For any type $T$ of order $d \ge 1$, we let  
	$\mathbf{e} = \mathbf{e}(T) = (e_{i_1}, e_{i_2}, \ldots, e_{i_t})$
be  the nonzero jumping indices $e_j$ of $T$, where $i_1 > i_2 > \cdots >i_t$. Notice that $\sum e_{i_j} = d$. We say that a $d \times d$ upper-triangular matrix $M=(m_{ij})$ {\em has shape $\mathbf{e}$} if the entries of the diagonal blocks of $M$ of sizes $e_{i_1} \times e_{i_1}, e_{i_2} \times e_{i_2}, \ldots,  e_{i_t} \times e_{i_t}$ are zero.
We let $\Mat_{\mathbf{e}}(\kk) \cong \AA^{n_{\mathbf{e}}}$ denote the affine space of $d \times d$ matrices of shape $\mathbf{e}$.
\end{definition}

\begin{example} \label{ex:shapee}
If $T=(1,2,3,4,5,3,3, 1,0, 0, \ldots)$, then degree is $d = 5$, the shape is $\mathbf{e}=(e_{i_1}, e_{i_2}, e_{i_3}) = (e_8,e_7,e_5) =(1,2,2)$, and $\Mat_{\mathbf{e}}(\kk) \cong \AA^{8}$ is the affine space of upper-triangular matrices of the form 
	$$
	\left(
	\scalemath{0.9}{\begin{array}{c|cc|ccc}
	0 & * & * & * & *  \\
	 \hline 0 & 0 & 0 & * & *   \\
	 0 & 0 & 0 & * & *  \\ 
	 \hline 0 & 0 & 0 & 0 & 0   \\
	 0 & 0 & 0 & 0 & 0  \\ 
	 \end{array}}\right).
	$$
\end{example}

A  dimension count shows that the matrices $\beta$ satisfying  (\ref{enum:betavanish})--(\ref{enum:beta0}) form an affine space of dimension $n_T := n - d - \sum_{i} \frac{e_i(e_i-1)}{2}$, whereas matrices of shape $\mathbf{e}$ with entries in $\kk$ form an affine space of dimension $n_{\mathbf{e}} := \sum_{j < k} e_{i_j} e_{i_k} = \frac{d^2 - \sum e_i^2}{2}$.

\begin{corollary}
\label{cor:2nd.Z_P}
For any normal pattern $P$ of type $T$ and order $d$, there is an isomorphism 
	\begin{equation}\label{eqn:2nd.Z_P}
\AA^{n_T}  \xrightarrow{\cong} Z_P, \quad \beta \mapsto I(\beta).
        \end{equation}
The minimal number of generators of $I(\beta)$ depends only on the image of $\beta$ under the coordinate projection, with fibers $\AA^{n_T - n_{\be}} = \AA^{n - \frac{d(d+1)}2}$, 
\[\AA^{n_T}\cong Z_P \to \AA^{n_{\mathbf{e}}} \cong \Mat_{\mathbf{e}}(\kk), \quad \beta \mapsto 
\overline{\beta(0)}:= \bigl(\beta_{ij}(0)\bigr)_{1 \leq i, j \leq d}
\]
obtained by removing the last row of zeros in $\beta(0)$: it is given by 
$d + 1 - \rank(\overline{\beta(0)})$. 
\end{corollary}
\begin{proof}
By virtue of Proposition \ref{prop:Iaparameters}, we only need to prove the claim about the minimal number of generators. Using the resolution \eqref{eq:Ibetaresolution} and Nakayama's Lemma, we see that it is given by $d+1 - \rank\left(\MM_P(0, 0) + \beta(0)\right)$.  As $\MM_P(0, 0) = 0$ and the rank of $\beta(0)$ is unchanged by removing the row of zeros at the bottom, the claim follows.
\end{proof}

\begin{example}
We consider an order $d=5$ stratum of $\Hilb_{22}(A)$ associated with the type $T = (1,2,3,4,5,3,3,1,0,0,\dots)$,  
where $A=\kk[[x,y]]$. The sequence of nonzero jumping indices is $\mathbf{e} = (e_8, e_7, e_5)=(1,2,2)$. The Young diagram $\Delta(P)$ and the values of $k_i$'s and nonzero $e_i$'s are illustrated in Figure \ref{figure:Delta(P)}.
\begin{figure}[htb]
    \centering
\begin{tikzpicture}[
BC/.style = {decorate, 
        decoration={calligraphic brace, amplitude=5pt, raise=1mm},
         very thick, pen colour={black}
            },
                    ]
\matrix (m) [matrix of math nodes,
             nodes={draw, minimum size=11mm, 
             anchor=center},
             column sep=-\pgflinewidth,
             row sep=-\pgflinewidth
             ]
{
|[draw=none]| y^5   & |[draw=none]|   & |[draw=none]|    & |[draw=none]|   & |[draw=none]|  & |[draw=none]|  & |[draw=none]|  & |[draw=none]| & |[draw=none]|& |[draw=none]|& |[draw=none]| & |[draw=none]| & |[draw=none]| k_5=0   \\
y^4  & |[draw=none]| x y^4   & |[draw=none]|    & |[draw=none]|   & |[draw=none]|  & |[draw=none]|  & |[draw=none]|  & |[draw=none]| & |[draw=none]|& |[draw=none]|& |[draw=none]| & |[draw=none]| & |[draw=none]| k_4=1    \\
y^3 & xy^3  &  |[draw=none]| x^2 y^3  & |[draw=none]|     & |[draw=none]|  & |[draw=none]|  & |[draw=none]|  & |[draw=none]| & |[draw=none]|& |[draw=none]|& |[draw=none]| & |[draw=none]| & |[draw=none]| k_3=2   \\
y^2  &  xy^2 & x^2y^2  & x^3y^2  &  x^4 y^2 &   |[draw=none]| x^5 y^2  & |[draw=none]|   & |[draw=none]|  & |[draw=none]| & |[draw=none]|& |[draw=none]|& |[draw=none]| & |[draw=none]| k_2=5\\
y  & xy &  x^2y & x^3y  &x^4y & x^5 y & |[draw=none]| x^6 y & |[draw=none]| & |[draw=none]| & |[draw=none]|& |[draw=none]|& |[draw=none]| & |[draw=none]| k_1=6   \\
1   & x & x^2 & x^3& x^4 & x^5 &x^6 &  x^7 &  |[draw=none]| x^8 & |[draw=none]| & |[draw=none]|& |[draw=none]| & |[draw=none]| k_0=8 \\
};
\draw[BC] (m-2-2.north) -- node[above right=2.2mm] {$e_5=2$} (m-3-3.east);
\draw[BC] (m-4-6.north) -- node[above right=2.2mm] {$e_7=2$} (m-5-7.east);
\draw[BC] (m-6-9.north) -- node[above right=2.2mm] {$e_8=1$} (m-6-9.east);
\end{tikzpicture}
    \caption{Young diagram $\Delta(P)$ in the case where $T = (1,2,3,4,5,3,3,1)$} along with the values of $k_i$'s and nonzero $e_i$'s.
    \label{figure:Delta(P)}
\end{figure}

The matrices $\beta$ satisfying  (\ref{enum:betavanish})--(\ref{enum:beta0}) form an affine space $\AA^{15}$, and for each such $\beta \in \AA^{15}$, the matrix $\MM_P + \beta$ of Proposition \ref{prop:Iaparameters} takes the form 
\begin{align*}
	\MM_P+ \beta
	 = 
	\left(
	\begin{array}{ccccc}
	-y + \beta_{11}^1x & \beta_{12}^0 & \beta_{13}^0+\beta_{13}^1x  +\beta_{13}^2x^2 & \beta_{14}^0 & \beta_{15}^0 \\
	 x^2 & -y & \beta_{23}^1x + \beta_{23}^2x^2 & \beta_{24}^0 & \beta_{25}^0 \\
	 0 & x & -y+\beta_{33}^1x + \beta_{33}^2 x^2 & \beta_{34}^0 & \beta_{35}^0 \\
	 0 & 0 & x^3 & -y & 0 \\
	 0 & 0 & 0 & x & -y \\
	 0 & 0 & 0 & 0 & x\end{array}
	 \right),
	\end{align*}
and the ideal $I(\beta)$ is generated by the $5\times 5$ minors of the above matrix; here, $\beta_{ij}^{k} \in \AA^1(\kk)$ are the coefficients of the polynomial $\beta_{ij} = \sum_k \beta_{ij}^k x^k$. The isomorphism \eqref {eqn:2nd.Z_P} in this case is the isomorphism $\beta \mapsto I(\beta)$, $\AA^{15} \xrightarrow{\cong} Z_{P}$. Moreover, the projection $(\beta \in \AA^{15}) \mapsto (\overline{\beta(0)} \in \AA^8)$ of Corollary \ref{cor:2nd.Z_P} takes the form:
$$
\beta = 	\left(\begin{array}{ccccc}
	\beta_{11}^1x & \beta_{12}^0 & \beta_{13}^0+\beta_{13}^1x  +\beta_{13}^2x^2 & \beta_{14}^0 & \beta_{15}^0 \\
	 0& 0 & \beta_{23}^1x + \beta_{23}^2x^2 & \beta_{24}^0 & \beta_{25}^0 \\
	 0 & 0 & \beta_{33}^1x + \beta_{33}^2 x^2 & \beta_{34}^0 & \beta_{35}^0 \\
	 0 & 0 & 0 & 0 & 0 \\
	 0 & 0 & 0 & 0 & 0 \\
	 0 & 0 & 0 & 0 & 0\end{array}
	 \right)
	\mapsto
	 \overline{\beta(0)}=
	\left(	\begin{array}{c|cc|cc}
	0 & \beta_{12}^0 & \beta_{13}^0 & \beta_{14}^0 & \beta_{15}^0 \\
	\hline
	 0 & 0 & 0 & \beta_{24}^0 & \beta_{25}^0 \\
	 0 & 0 &0 & \beta_{34}^0 & \beta_{35}^0 \\
	 \hline
	 0 & 0 & 0 & 0 & 0 \\
	 0 & 0 & 0 & 0 & 0
	 \end{array}
	 \right).
	 $$
\end{example}

\def\ua{\underline{a}}

In light of Corollary \ref{cor:2nd.Z_P}, we need to describe the degeneracy locus of matrices of shape $\be$ and of given rank. We first consider a slightly more general problem. Given a non-decreasing function $\Gamma \colon \{1, 2, \dots, d\} \to \{0, 1, 2, \dots, d\} $, we say that a $d \times d$-matrix $M$ is of type $\Gamma$ if $M_{ij} = 0 $ for $i > \Gamma(j)$.  

\begin{remark} \label{rem:shape-vs-type}
Let $\be$ be a collection of jumping indices, and let $\Gamma$ be the function with
$\Gamma(i) = e_1 + \dots + e_k$ for $e_1 + \dots + e_ k < i \leq e_1 + \dots + e_{k+1}$ and $0 \le k \le s-1$. 
Then a matrix of shape $\be$ if and only it is of type $\Gamma$. (So in Example~\ref{ex:shapee} we would set $\Gamma(1) = 0$, $\Gamma(2) = \Gamma(3) = 1$ and $\Gamma(4) = \Gamma(5) = 3$.)
But functions $\Gamma$ not of this form do not correspond to a shape $\be$.
\end{remark}

\begin{lemma} \label{lem:Gammamatrixdegeneracy} Let $\DD_R^\Gamma$ be the locus of $d \times d$-matrices of type $\Gamma$ and rank $R$. Then $\DD_R^\Gamma$ is non-empty if and only if, for all $1 \le k \le d$, we have
\begin{equation} \label{eq:DrG-nonempty}
    \Gamma(k) -k \ge R - d.
\end{equation}\

For each sequence $1 \le a_1 < a_2 < \dots < a_R \leq d$, denoted by $\ua$, let
\begin{equation}\label{eq:DrG-dimension}
\rho^\Gamma(\ua) := Rd - \frac{R(R-1)}2 + \sum_{i=1}^R \left(\Gamma(a_i) - a_i\right).
\end{equation}
If non-empty,  the dimension of $\DD_R^\Gamma$  is the maximum of $\rho^\Gamma(\ua)$ for all sequences $\ua$ as above satisfying $\Gamma(a_i) \ge i$ for all $1 \leq i \leq R$.
\end{lemma}
\begin{proof}
Given a matrix $M$ of rank $R$, we let $\ua = (a_1, \dots, a_R)$ be the sequence describing its row-echelon form:  $a_i$ is the minimal number such that the first $a_i$ columns of $M$ have rank $i$. The maximal possible such sequence is given by $a_i = d+i - R$, e.g.~for the matrix $N_R$ obtained as the $R\times R$-unit matrix a the top right, filled with zeros in the rows below and columns to the left; every other possible sequence satisfies $a_i \le d+i - R$.

If $M$ is of type $\Gamma$, then the first $k$ columns have rank at most $\Gamma(k)$, as they only contain zeros below the row $\Gamma(k)$. Therefore, $\Gamma(a_i) \ge i$ for all $i$. As $\Gamma$ is monotone, combined with the previous inequality we get
\[
\Gamma(d + i - R) \ge \Gamma(a_i) \ge i,
\]
for all $i = 1, \dots , R$, i.e.~inequality \eqref{eq:DrG-nonempty} for $k > d - R$. As $\eqref{eq:DrG-nonempty}$ is trivially satisfied for $k \le d-R$ due to $\Gamma(k) \ge 0$, this shows that $\eqref{eq:DrG-nonempty}$ is necessary for $\DD^\Gamma_R$ to be non-empty. On the other hand, the matrix $N_R$ shows that it is also sufficient.

Given $M \in \DD_R^{\Gamma,\ua}$,  there is a basis $v_1, v_2, \dots, v_R$ (the rows in a row-echelon form of $M$)  of the row space of $M$ such that the first non-zero entry of $v_i$ is at position $a_i$; this basis is well-defined up to the action of upper triangular $R \times R$-matrices, and the dimension of the space of such bases is given by $\sum_i (d + 1 - a_i)$. Every row of $M$ is a linear combination of $v_1, \dots, v_R$, where $v_i$ can contribute to the first $\Gamma(a_i)$ rows; moreover, a generic choice of such a linear combination yields a matrix of rank $R$ and type $\Gamma$.  Thus $\DD_R^{\Gamma,\ua}$ is irreducible and of dimension
\[
\sum_{i=1}^R (d+1 - a_i) - \frac{R(R+1)}2 + \sum_{i=1}^R  \Gamma(a_i) = \rho^\Gamma(\ua).
\]
As $\DD_R^\Gamma$ is the union of $\DD_R^{\Gamma,\ua}$ for all possible sequences $\ua$, this proves the claim.
\end{proof}

\begin{corollary}
\label{cor:deg.loci:Mat_e}
Let $d \ge 2$ be an integer, $\mathbf{e}= (e_1, \ldots, e_s)$ be an ordered partition of $d$ of length $s = {\rm length}(\be)$.  For each $R \ge 0$, let $\DD_R(\Mat_{\mathbf{e}}(\kk)) = \{M \in \Mat_{\mathbf{e}}(\kk) \mid \rank M = R\}$ denote the degeneracy locus of matrices of shape $\be$ and rank $R$. Then 
\begin{enumerate}
	\item \label{enum:maternonempty} $\DD_R(\Mat_{\mathbf{e}}(\kk))\neq \emptyset$ if and only if $ R \le d - \max\{e_j\}$.
        \item \label{enum:mater-dimension}
        The dimension is bound by
        \begin{equation} \label{eq:DrMatedim}
            \dim \DD_R(\Mat_{\mathbf{e}}(\kk)) \le R \left(\frac{2d-R-1}{2}\right).
        \end{equation}
        \item \label{enum:mater-maxdim}
        If ${\rm length}(\mathbf{e}) \ge R+1$, then \eqref{eq:DrMatedim} is an equality.
\end{enumerate}
\end{corollary}
\begin{proof}
We  apply Lemma \ref{lem:Gammamatrixdegeneracy} with $\Gamma$ as defined in Remark \ref{rem:shape-vs-type}.
The minimal possible value of $\Gamma(k)-k$ is given by $-\max\{e_j\}$; substituting this value in equation \eqref{eq:DrG-nonempty} proves (\ref{enum:maternonempty}).

Since $\Gamma(k) \leq k-1$ for all $1 \le k \le d$, the last term in formula \eqref{eq:DrG-dimension} satisfies 
\[ \sum_{i = 1}^R \Gamma(a_i) - a_i \le -R,\]
which proves (\ref{enum:mater-dimension}). Moreover, equality $\Gamma(k) = k-1$ holds only for $k \in \{e_1 + 1, e_1 + e_2 + 1, \dots, e_1 + \dots + e_{s-1} + 1 \}$; thus equality in \eqref{eq:DrMatedim}  is obtained if  and only if the sequence $\ua$ is a subsequence of $(e_1 + 1, e_1 + e_2 + 1, \dots, e_1 + \dots + e_{s-1} + 1)$, which is possible if ${\rm length}(\be) \ge R+1$. This proves (\ref{enum:mater-maxdim}).
\end{proof}

\begin{remark} \label{rem:DrMate-Grassmann}
    We note the special case of ${\rm length}(\be) = 2$, which includes the case of the Grassmannian stratum of Example \ref{ex:grassmannstratum}. 
    
    In this case, matrices in $\Mat_{\be}(\kk)$ are just $e_1 \times e_2$-matrices, extended by rows of zeros and columns on the bottom and left to form a $d \times d$-matrix.  Thus $\DD_r(e_1, e_2)$ is non-empty if and only if $r \le \min\{e_1, e_2\}$, in which case it has dimension $r(d-r)$. 
\end{remark}

\subsection{Local Brill--Noether theory}
\label{sec:localBN}

The following is our main result regarding the Brill--Noether loci on each Hilbert--Samuel stratum $Z_T$. The smallest possible minimal number of generators for ideals in $Z_P$ had previously been described in \cite[Theorem 4.3]{Ia} and \cite[Proposition III.2.1]{Bri}.

\begin{theorem}[Brill--Noether for Hilbert--Samuel Strata] 
\label{thm:local}
Let $T$ be a type of order $d$ and $|T|=n$, and let $Z_T \subset \Hilb_n(\kk[[x,y]])$ denote the associated stratum. Assume that either ${\rm char}(\kk) =0$ or ${\rm char}(\kk) \ge |T|=n$. Let $r_{\rm min}: = \max_j\{e_j\}$, where $e_j$ are the jumping indices of $T$. We let $\mathbf{e} = \mathbf{e}(T)$ be the list of nonzero jumping indices $e_j$ of $T$ as in Definition \ref{def:Mat_e}, and let ${\rm length}(\mathbf{e})$ denote the length of the sequence $\mathbf{e}$. For each $r>0$, we consider the Brill--Noether locus of the stratum $Z_T$,
	$$\BN_{=r} (Z_T) : = \{I \in Z_T \mid \dim_{\kk}(I \otimes_{A} \kk) = r +1\} \subset Z_T,$$
where $\dim_{\kk}(I \otimes_{A} \kk)$ is the minimal number of generators of an ideal $I$. Then:
\begin{enumerate}
	\item 
        \label{enum:BNrZT-nonempty}
      We have  $\BN_{=r} (Z_T) \ne \emptyset$ if and only if $r_{\rm min} \le r \le d$. 
	\item 
       \label{enum:BNrZT-dimbound}
        For all  $r_{\rm min} \le r \le d$, we have
		\begin{equation} \label{eq:dimBnrZT}
                \dim \BN_{=r} (Z_T) \le n - \frac{r(r+1)}{2} - (d- r).
		\end{equation} 
        \item \label{enum:BNrZT-dimequality}
         If 
        ${\rm length}(\mathbf{e}) \ge (d-r)+1$, then $\BN_{=r}(Z_T)$ is nonempty, and \eqref{eq:dimBnrZT} is an equality.
	\item
        \label{enum:BNd-ZT} In particular, 
        $\BN_{=d} (Z_T)$ is nonempty and of dimension $n - \frac{d(d+1)}{2}$.
\end{enumerate}
 \end{theorem}

\begin{proof}
By virtue of Proposition \ref{prop:affinecover} we only need to prove the  results on the affine chart $Z_P$. 
Corollary \ref{cor:2nd.Z_P} implies that the isomorphism \eqref{eqn:2nd.Z_P} induces an isomorphism
    $$\DD_{d-r}(\Mat_{\mathbf{e}}(\kk)) \times \AA^{n-\frac{d(d+1)}{2}} \xrightarrow{\cong}\BN_{=r}(Z_P) \quad \text{for all} \quad r \ge 0,$$
where $\DD_{d-r}(\Mat_{\mathbf{e}}(\kk))$ is as defined in Corollary \ref{cor:deg.loci:Mat_e}, and
$\BN_{=r}(Z_P) = \BN_{=r}(Z_T) \cap Z_P$.

Then the claims (\ref{enum:BNrZT-nonempty}), (\ref{enum:BNrZT-dimbound}) and (\ref{enum:BNrZT-dimequality}) follow directly from the corresponding claims in Corollary \ref{cor:deg.loci:Mat_e} after substituting $R = d-r$. Finally, (\ref{enum:BNd-ZT}) is a special case of (\ref{enum:BNrZT-dimequality}).
\end{proof}

\begin{example}
      Consider the Grassmannian stratum in $\Hilb_{\frac{d(d+1)}2 + \ell}(S)$ of Example \ref{ex:grassmannstratum} for given $d$ and $\ell$, and let $\Gr$ be the corresponding type. In light of Remark~\ref{rem:DrMate-Grassmann}, we can be more precise about the loci $\BN_r(Z_\Gr)$. The locus $\BN_r(Z_\Gr)$ is non-empty if and only if $\max\{\ell, d-\ell\} \le r \le d$; in this case $\BN_r(Z_\Gr)$ is irreducible and of dimension $\ell + r(d-r)$. In particular, $\BN_d(Z_\Gr)$ is non-empty and has dimension $d$.
\end{example}

From Theorem \ref{thm:local}, we obtain the following:

\begin{corollary}\label{cor:local} 
Let $A = \kk[[x,y]]$, where $\kk$ is a field, and let $n \ge 2$ be an integer. Assume that ${\rm char} (\kk)=0$ or ${\rm char}(\kk) \ge n$. Let $\Hilb_n(A)$ denote the Hilbert scheme of $n$ points on $R=\kk[[x,y]]$ with the reduced scheme structure. For any integer $r \ge 0$, we define 
    $$\rho^{\rm loc}_{r,n}: = n - \frac{r(r+1)}2.$$
We let $\BN_{r,n}^{\rm loc}$ denote the Brill--Noether locus $\BN_{r}(\Hilb_n(A)) = \{ I \in \Hilb_n(A) \mid \dim_{\kk}(I \otimes_{A} \kk)  \ge r +1 \}$. Then:
\begin{enumerate}
	\item \label{cor:local-1} $\BN_{r,n}^{\rm loc} \neq \emptyset$ if and only if $\rho_{r,n}^{\rm loc} \ge 0$.
        \item \label{cor:local-2} If the conditions of \eqref{cor:local-1} are satisfied, then $\dim \BN_{r,n}^{\rm loc} = \rho_{r,n}^{\rm loc}$. 
        \item \label{cor:local-3}
        If $\rho_{r,n}^{\rm loc} = 0$, then $\BN_{r,n}^{\rm loc} \cong \{\fom^r\}$ is a point.
\end{enumerate}
\end{corollary}

\begin{proof} 
By virtue of Theorem \ref{thm:stratification}, $\BN_{r,n}^{\rm loc}$ is the union of $\BN_{r,n}(Z_T)$, where $T$ runs through all possible types with $|T|=n$.
The condition $\rho_{r,n}^{\rm loc} \ge 0$ is equivalent to the condition that there is a type $T$ with $|T|=n$ and order $r$. For such a type $T$, Theorem \ref{thm:local} \eqref{enum:BNd-ZT} implies that $\BN_r(Z_T)$ is nonempty and has dimension $\rho_{r,n}^{\rm loc}$.

Conversely, if $\BN_{r,n}^{\rm loc} \neq \emptyset$, then there exists a type $T$ such that $\BN_{r}(Z_T) \neq \emptyset$. By Theorem \ref{thm:local} \eqref{enum:BNrZT-nonempty}, such a type $T$ has order $\ge r$, which implies $\rho_{r,n}^{\rm loc} \ge 0$. This proves claim \eqref{cor:local-1}.

Now assume  the conditions of claim \eqref{cor:local-1}. If a type $T$ has order $< r$, then $\BN_r(Z_T)$ is empty by Theorem \ref{thm:local}.\eqref{enum:BNrZT-nonempty}; if $T$ has order $r$, then $\dim \BN_{r}(Z_T) = \rho_{r,n}^{\rm loc}$ by Theorem \ref{thm:local}.\eqref{enum:BNrZT-dimequality}; and if a type $T$ has order $>r$, then $\dim \BN_r(Z_T)$ has dimension strictly smaller than $\rho_{r,n}^{\rm loc}$ by Theorem \ref{thm:local}.\eqref{enum:BNrZT-dimbound}. As there always exists a type of order $r$, this proves claims \eqref{cor:local-2}. And if $\rho_{r, n}^{\rm loc}= 0$, then the only possible type of order $r$ is $(1, 2, \dots, r, 0, 0, \dots)$, corresponding to the unique ideal $\fom^r$; this proves (\ref{cor:local-3}).
\end{proof}

\begin{remark} \label{rem:VeronoeseP1}
        When $\rho_{r,n}^{\rm loc}=1$, $\Hilb_n(A)$ consists of precisely one degree $r$ stratum $Z_T$ which parametrizes ideals of type $T= (1,2,\dots,r,1,0,0 \dots)$, i.e., $Z_T \cong \PP^r$ is the Grassmannian stratum of Example \ref{ex:grassmannstratum} in the case where $\ell=1$. In this case, $\BN_{r,n}^{\rm loc} = \BN_{r}(Z_T) \cong \PP^1$ is the image of $r$th Veronese embedding $\nu_r \colon \PP^1 \hookrightarrow \PP^r \cong Z_T$. 
        To see this, let $P = P(x,y)$ be a normal pattern associated with $T$. Then ideals of $Z_P$ are given by $I = (f_0, f_1, \ldots, f_r)$, where $f_0 = x^{r+1}$, $f_{1} = x^{r-1} y - a_{1} x^r$, $\ldots$, $f_{r} =  y^r - a_{r} x^r$, $(a_1,\ldots,a_r) \in \AA^r$. Each such ideal $I$ has a presentation $R^{r} \xrightarrow{M} R^{r+1} \to I$, where 
	$$
	M = \left(\begin{array}{ccccc}
	-y+a_{1}x & a_{2}-a_{1}^2 & a_{3}-a_{1}a_{2} & \cdots & a_{r}-a_{1}a_{r-1}  \\
	x^2 & -y-a_{1}x & -a_{2} x& \cdots & -a_{r-1} x \\
	 & x & -y & 0 & \cdots \\  &   & x & \cdots & 0 \\  &   &   & \cdots & -y \\  &   &   &   & x\end{array}\right).
	$$
 Therefore, $I\in \BN_{r}(Z_P)$ if and only if $M|_{(0,0)}=0$, if and only if
    \[a_2-a_1^2 = a_3 - a_1 a_2 = \cdots = a_r - a_1a_{r-1}  = 0.\]
Applying the action of $\GL_2(\kk)$, we see that $\BN_r(Z_T) \subset \PP^r$ is smooth, closed, and one-dimensional, and is the closure of the curve $\{(a_1,a_1^2, \cdots, a_1^r) \mid a_1 \in \AA^1\} \subset \AA^r$ in $\PP^r$.
\end{remark}

\subsection{Global Brill--Noether theory}
\label{sec:globalBN}
Before applying our local results in the global setting, we note that our Brill-Noether loci have a well-defined expected dimension:
\begin{lemma} \label{lem:expecteddim}
    The locus $\BN_{r, n}$ is everywhere of dimension at least
    \[ \rho_{r, n} = 2n + 2 - r(r+1).
\]
If this expected dimension is achieved at each point, then $\BN_{r,n}$ is Cohen--Macaulay.
\end{lemma}
\begin{proof}
Since the universal subscheme $Z_n \subset \Hilb_n(S) \times S$ is finite and flat over the smooth scheme $\Hilb_n(S)$, it is Cohen--Macaulay. Let $z=(I,p) \in \Hilb_n(S) \times S$. By the Hilbert–Burch theorem (\cite[Theorem 20.15]{Ei}), the ideal $I_{Z_n}$ admits a resolution of length two of the form
    \[\sO_{\Hilb_n(S) \times S, z}^{k} \xrightarrow{M} \sO_{\Hilb_n(S) \times S,z}^{k+1} \to I_{Z_n,z}\]
near $z$, where $\ell \ge 1$ is an integer.

Since $I_{Z_n,z} \otimes_{\sO_{\Hilb_n \times S,z}} \kappa(z) \cong I_p \otimes_{\sO_{S,p}} \kappa(p)$, the Brill--Noether locus $\BN_{r,n}$ at $z=(I,p)$ coincides with the degeneracy locus where the $k \times (k+1)$ matrix $M$ has rank $\le k-r$.
Therefore, it is a closed subset of codimension at most $r(r+1)$. Moreover, when the maximal codimension $r(r+1)$ is achieved, the Brill--Noether locus $\BN_{r,n}$ is Cohen--Macaulay at the point $z$ (see, for example, \cite[Exercise 10.9 \& Theorem 18.18]{Ei}).
\end{proof}

\begin{proof}[Proof of Theorem \ref{thm:main} under the condition ${\rm char}(\kk) = 0$ or ${\rm char}(\kk) \ge n$]Consider the stratification of $\BN_{r,n}$ by the multiplicity $m_p(I) = \dim_{\kk}(\sO_{S,p}/I_{S,p})$ of $I$ at $p$:
	$$\BN_{r,n} = \bigsqcup_{m=1}^n \BN_{r,n}^{(m)}, \qquad  \BN_{r,n}^{(m)} = \{(I,p) \in \BN_{r,n} \mid m_p(I) = m\}.$$
By semicontinuity, each $\BN_{r,n}^{(m)}$ is a locally closed subset.

Consider the projection $\BN_{r, n}^{(m)} \to S$ that sends $(I, p)$ to $p$. Choosing local parameters $x, y$ at a point $p \in S$ identifies the fiber of this projection over $p$ with $\BN_{r,m}^{\rm loc} \times \Hilb_{n-m}(S \backslash \{p\})$. 

From Corollary \ref{cor:local} \eqref{cor:local-1}, we obtain that $\BN_{r,n}^{(m)} \ne \emptyset$ if and only if $\frac{r(r+1)}{2} \le m \le n$. Consequently, $\BN_{r,n} \ne \emptyset$ if and only if $\frac{r(r+1)}{2} \le n$, if and only if $\rho_{r,n} \ge 2$. If $\BN_{r,n}^{(m)} \ne \emptyset$, by Corollary \ref{cor:local} \eqref{cor:local-2}, we have 
	$$\dim \BN_{r,n}^{(m)}  = \dim \BN_{r,m}^{\rm loc} + 2(n-m+1) = 2n +2  - m  -\frac{r(r+1)}{2}.$$
So $\dim \BN_{r,n}^{(m)}$ is strictly decreasing with respect to the variable $m \in [\frac{r(r+1)}{2}, n]$, and achieves its maximum if and only if $m = m_{\min} := \frac{r(r+1)}{2}$. In this case, 
\[ \dim \BN_{r,n}^{(m_{\min})}=2n+2 - r(r+1) = \rho_{r,n},
\]
and  the fiber of $\BN_{r,n}^{(m_{\min})}$ under the projection to $S$ is $\BN_{r,m_{\min}}^{\rm loc} \times \Hilb_{n-m_{\min}}(S \backslash \{p\}) = \{\mathrm{point}\} \times \Hilb_{n-m_{\min}}(S \backslash \{p\})$ by Corollary~\ref{cor:local}.\eqref{cor:local-3}. In particular,   $\BN_{r,n}^{(m_{\min})}$ is irreducible.

By Lemma \ref{lem:expecteddim}, this implies that $\BN_{r,n}^{(m_{\min})}$ of $\BN_{r,n}$ is irreducible of the expected dimension, as every lower-dimensional stratum is contained in the closure of this one.
\end{proof}

\section{Birational correspondences between Brill-Noether loci}
\label{sec:biration}

In this section, we give a second proof of Theorem~\ref{thm:main} using birational correspondences between different Brill-Noether loci defined by certain nested Hilbert schemes. For simplicity, we continue to assume that $\kk$ is algebraically closed.

Recall from the introduction that
\[ 
\BN_{r, n} = \left\{ (I,p) \mid \dim_\kk (I \otimes_{\sO_S} \kappa(p)) \ge r+1\right
\} \subset  \Hilb_n(S) \times S,
\]
 for $n \ge 1, r \ge 0$, and let $\BN_{=r, n}$ be the open subset
\[ 
\BN_{=r, n} = \left\{ (I,p) \mid \dim_\kk (I \otimes_{\sO_S} \kappa(p)) = r+1\right
\} \subset  \Hilb_n(S) \times S.
\]
Recall from Lemma~\ref{lem:expecteddim} (the only result from Section~\ref{sec:local} that we will use) that $\BN_{r,n}$ has expected 
 dimension $\rho_{r, n}$. 
We will prove by induction on $n$ that $\BN_{r,n}$ and $\BN_{=r,n}$ have expected dimension 
and are non-empty if and only if $\rho_{r, n} \ge 2$. Note that
$\rho_{r-1, n} = \rho_{r, n+r}$, which is the first hint of a relation between $\BN_{r-1, n}$ and $\BN_{r, n+r}$. 
The key role in this relation, and in our second proof of Theorem~\ref{thm:main}, is played by the following  nested Hilbert scheme:

\begin{definition}
For $n \ge 1, 1\le r \le n$ we let $\Hilb_{n-r,n}^{\dagger}(S)$ be the $r$-step nested Hilbert scheme 
\[
\Hilb_{n-r,n}^{\dagger}(S) : = \{ I_{n} \subset I_{n-r}  \mid
\exists p \in S , I_{n-r}/I_{n} \cong \kappa(p)^{\oplus r} \} \subset \Hilb_{n-r}(S) \times \Hilb_{n}(S).
\]
We let $\pi_1, \pi_2$ be the natural projections:	\begin{equation}\label{eqn:nested.Hilbert}
	\begin{tikzcd}[row sep= 2.2 em, column sep={8 em,between origins}]
		 & \Hilb_{n-r,n}^{\dagger}(S)  \ar{dr}{\pi_2} \ar{dl}[swap]{\pi_1} & \\
		 \Hilb_{n-r}(S) \times S & &  \Hilb_{n}(S) \times S
	\end{tikzcd}
	\end{equation}
\end{definition}

It was proven in \cite[Lemma 5.6]{J20} that both $\pi_1$ and $\pi_2$ can be described as relative Grassmannian bundles, for the ideal $\sI_{Z}$ of the universal subscheme $Z \subset \Hilb_{n-r}(S) \times S$ in case of $\pi_1$, and for $\omega_{Z}$ in case of $\pi_2$. 
For our purposes, we only need the following consequences of 
\cite[Lemma 5.6]{J20}, for which we give a self-contained proof for convenience:
\begin{lemma} \label{lem:Grassmannian-fibers}
\begin{enumerate}
    \item \label{enum:pi1image} The image of $\pi_1$ is $\BN_{r-1,n-r}$. 
    \item \label{enum:pi1fibers} Moreover, if $I_n \in \BN_{=r'-1, n-r}$ for some $r' \ge r$, then the fiber $\pi_1^{-1}(I_n)$ is isomorphic to the Grassmannian $\Grass(r, r')$. In particular, it is an isomorphism over $\BN_{=r-1, n-r}$.
    \item  \label{enum:pi2image} The image of $\pi_2$ is $\BN_{r, n}$, and $\pi_2$ is an isomorphism over $\BN_{=r, n}$.  If $r > n$, then $\BN_{r, n}$ is empty.
    \end{enumerate}
\end{lemma}
\begin{proof}
By definition, closed points of $\Hilb_{n-r,n}^{\dagger}(S)$ are in 1:1-correspondence with short exact sequences of the form
\begin{equation} \label{eq:nested-ses}
0 \to I_{n} \to I_{n-r} \to \kappa(p)^{\oplus r} \to 0.
\end{equation}
First consider $I_{n-r}$ fixed. Then 
such sequences correspond to surjections
$I_{n-r} \onto \kappa(p)^{\oplus r}$ up to the automorphism $\GL_r(\kappa(p))$ of $\kappa(p)^{\oplus r}$, which in turn are given by surjective maps
$I_{n-r} \otimes_{\sO_S} \kappa(p) \onto \kappa(p)^{r}$ of $\kk \cong \kappa(p)$-vector spaces; this proves both (\ref{enum:pi1image}) and (\ref{enum:pi1fibers}).

Now we apply $\Hom(\blank, \kappa(p))$ to the short exact sequence \eqref{eq:nested-ses}. As $\Ext^2(I_{n-r}, \kappa(p)) = 0$, we obtain a surjection
\[
 \Ext^1(I_{n}, \kappa(p)) \onto \Ext^2\left(\kappa(p)^{\oplus r}, \kappa(p)\right) \cong \kappa(p)^r. 
\]
Since $\chi(I_{n}, \kappa(p)) = 1$ and all higher Ext vanish by Serre duality, this gives
\[ \dim_{\kappa(p)} (I_n \otimes_{\sO_S} \kappa(p)) = \dim_{\kappa(p)} \Hom(I_n, \kappa(p))
= \dim_{\kappa(p)} \Ext^1(I_n, \kappa(p))+1 \ge r+1.
\]
This shows that the image of $\pi_2$ is contained in $\BN_{r, n}$.

Conversely, if $I_{n} \in \BN_{r, n}$, then $\dim_{\kappa(p)} \Ext^1(I_{n}, \kappa(p)) = \dim_{\kappa(p)} \Ext^1(\kappa(p), I_{n}) \ge r$ by the same reasoning. Any $r$-dimensional subspace of $\Ext^1(\kappa(p), I_{n})$ defines an extension
\[
0 \to I_{n} \to J \to \kappa(p)^{\oplus r} \to 0
\]
with $\Hom(\kappa(p), J) = 0$. Thus $J$ is a torsion-free sheaf of rank one and trivial determinant on a smooth surface, and therefore an ideal sheaf. This shows that $r \le n$, and that $I_n$ is contained in the image of $\pi_2$. Moreover, if $I_{n} \in \BN_{=r, n}$, this short exact sequence is unique. This  concludes the proof of (\ref{enum:pi2image}). 
\end{proof}

\begin{theorem} \label{thm:birational}  \begin{enumerate}
    \item \label{enum:BNrn-nonempty}
For each $n \ge 0$ and $r \ge 0$, the locus $\BN_{r, n}$ is non-empty if and only if its expected dimension satisfies $\rho_{r, n} \ge 2$.
\item \label{enum:BNrn-expdim}
If this condition is satisfied, then $\BN_{r, n}$ is irreducible, of expected dimension $\rho_{r, n}$, and birational to
$\Hilb_i(S) \times S$, for $i = \frac{\rho_{r, n}}2 - 1$. 
\item \label{enum:birational}
Moreover, if $n \ge 1$, $1 \le r \le n$ and $\rho_{r, n}\ge 2$, then the diagram \eqref{eqn:nested.Hilbert} induces a diagram
	\begin{equation} \label{eq:bir-correspondence}
	\begin{tikzcd}[column sep = 3 em, row sep = 1.5 em]
		 & \Hilb_{n-r,n}^{\dagger}(S)  \ar{dr}{\overline{\pi}_2} \ar{dl}[swap]{\overline{\pi}_1} & \\
		 \BN_{r-1,n-r} & &  \BN_{r,n},
	\end{tikzcd}
	\end{equation}
where $\overline{\pi}_1$ and $\overline{\pi}_2$ are isomorphisms over the open and dense loci $\BN_{=r-1, n-r}$ and $\BN_{=r, n}$.
\end{enumerate}
\end{theorem}

\begin{proof}
We will give a proof by induction on $n$.
If $n=0$,  $\BN_{0, 0} = \Hilb_0(S) \times S$ and $\BN_{r, 0} = \emptyset$ for $r \ge 1$, as claimed.

For the induction step, by Lemma \ref{lem:Grassmannian-fibers}.(\ref{enum:pi2image}), we can assume $r \le n$.  For $r = 0$, we have $\BN_{0,n} = \Hilb_n(S) \times S$, matching our claims. In the remaining cases $1 \le r \le n$ we use diagram \eqref{eqn:nested.Hilbert} and apply Lemma \ref{lem:Grassmannian-fibers}. Consider the map $\overline{\pi}_1 \colon \Hilb_{n-r, n}^\dagger \to \BN_{r-1, n-r}$ and the stratification
\[ \BN_{r-1, n-r} = \bigcup_{r' \ge r} \BN_{=r'-1, n-r}  \]
of the image.
By induction assumption, $\BN_{=r'-1, n-r}$ is irreducible and of expected dimension $ \rho_{r'-1, n-r} = 2n - 2r + 2 - (r'-1)r'$, and non-empty iff $ \rho_{r'-1, n-r} \ge 2$. In particular, $ \Hilb_{n-r, n}^\dagger$, and thus $\BN_{r, n}$, is non-empty if and only if $\rho_{r,n} = \rho_{r-1, n-r} \ge 2$, proving (\ref{enum:BNrn-nonempty}).

By Lemma \ref{lem:Grassmannian-fibers}.(\ref{enum:pi1fibers}), the preimage  $\overline{\pi}_1^{-1}(\BN_{=r'-1, n-r})$ has dimension
\begin{align*}
2n - 2r + 2 - (r'-1)r' + r(r' - r)
&= 2n + 2 - r(r+1) - (r'-1)(r'-r)  \\
&= \rho_{r, n} -  (r'-1)(r'-r)   \leq \rho_{r, n},
\end{align*}
with equality only if $r' = r$. As $\BN_{r-1, n-r}$ is irreducible by the induction assumption, $\Hilb_{n-r,n}^{\dagger}(S)$ has one irreducible component of dimension $\rho_{r, n}$, birational to $\BN_{r-1, n-r}$; should $\Hilb_{n-r,n}^{\dagger}(S)$ be reducible, every other component has smaller dimension. Since $\BN_{r, n}$ is the image of $\Hilb_{n-r,n}^{\dagger}(S)$ by Lemma \ref{lem:Grassmannian-fibers}.(\ref{enum:pi2image}), and has dimension at least $\rho_{r, n}$ at every point by Lemma \ref{lem:expecteddim}, it  is irreducible, of expected dimension, and birational to $\BN_{r-1, n-r}$.  This proves both claim (\ref{enum:BNrn-expdim}) and, when combined with Lemma \ref{lem:Grassmannian-fibers}.(\ref{enum:pi2image}), claim (\ref{enum:birational}).
\end{proof}

\begin{remark} \label{rem:nestedHilbert}
In the situation of Theorem \ref{thm:birational}.\eqref{enum:BNrn-expdim}, $\Hilb_{n-r,n}^{\dagger}(S)$ is also irreducible, of expected dimension $\rho_{r, n}$, and birational to $\Hilb_i(S) \times S$.  Since $\sI_{Z_{n-r}}$ locally admits a resolution of the form $\sO_{\Hilb_{n-r}(S) \times S}^k \xrightarrow{M} \sO_{\Hilb_{n-r}(S) \times S}^{k+1}$ for some $k \ge 1$ (see the proof of Lemma \ref{lem:expecteddim}), we can locally identify points of $\Hilb_{n-r,n}^{\dagger}(S)$ with surjections $\shO_{\Hilb_{n-r}(S) \times S}^{k+1} \to \kappa(z)^r$ that become zero when composed with $M$. This locally identifies $\pi_1$ with the zero locus of the section of an $rk$-dimensional vector bundle on a relative $\Grass(r, k+1)$-bundle over $\Hilb_{n-r}(S) \times S$ (see e.g.~\cite[Proposition 4.19]{J22b} for a general description of this phenomenon). Thus it has dimension at least  $2(n-r)+2+r(k+1-r) - rk = \rho_{r,n}$ at every point, and so the irreducible component described in the proof of Theorem~\ref{thm:birational} is the only one. (See \cite[Lemma 7.16]{J22a} for a similar argument, for projectivisations rather than relative Grassmannians.)
\end{remark}

\begin{example}
    When $\BN_{r, n-r}$ and $\BN_{r+1, n}$ are empty, then by Lemma \ref{lem:Grassmannian-fibers} both $\overline{\pi}_1$ and $\overline{\pi}_2$ in \eqref{eq:bir-correspondence} are isomorphisms, and so $\BN_{r-1, n-r} \cong \BN_{r, n}$.
    For any fixed expected dimension  $\rho_{r, n}$, this holds as soon as $2r \ge \rho_{r, n}$:
    \begin{enumerate}
        \item When the expected dimension is two, the locus $\BN_{1, 1} \subset \Hilb_1(S) \times S \cong S \times S$ is the diagonal $\Delta_S$, and $\BN_{1, 1} \cong \BN_{2, 3} \cong \dots \cong \BN_{d, \frac{d(d+1)}2} \cong S$ for all $d$. Concretely, $ \BN_{d, \frac{d(d+1)}2}$ is the locus of pairs $(\fom_p^d, p)$.
        \item When the expected dimension is four, we have $\BN_{0,1 } =  \Hilb_1(S) \times S \cong S \times S$, whereas $\BN_{1, 2} \subset \Hilb_2(S) \times S$ is the universal subscheme, isomorphic to $\Bl_{\Delta_S}(S \times S)$, and $\BN_{1, 2} \cong \BN_{2, 4} \cong \dots \cong \BN_{d, \frac{d(d+1)}2 +1 }$ for $d \ge 1$.
        \item When the expected dimensions is six, we have that $\BN_{0, 2} \cong \Hilb_2(S) \times S$ and 
        $\BN_{1, 3} \subset \Hilb_3(S) \times S$ is the universal subscheme. In the diagram \eqref{eq:bir-correspondence} for $\Hilb_{2, 3}^\dagger(S)$, both morphisms $\overline{\pi}_1$ and  $\overline{\pi}_2$ are non-trivial. In the next step, 
        $\overline{\pi_1} \colon \Hilb_{3, 5}^\dagger(S) \to \BN_{1, 3}$ is a $\Gr(2,3) = \PP^2$-bundle over $\BN_{2,3} \cong S$, and an isomorphism on the complement, whereas $\overline{\pi}_2 \colon \Hilb_{3, 5}^\dagger(S) \to \BN_{2, 5}$ is an isomorphism, with $\BN_{2, 5} \cong \BN_{3, 8} \cong \dots \cong \BN_{d,\frac{d(d+1)}2 + 2 }$ for $d \ge 3$.
    \end{enumerate}
\end{example}


\end{document}